\documentclass[reqno]{amsproc}

\usepackage{amstext,amsmath,amssymb,amsfonts}
\usepackage[latin1]{inputenc}
\usepackage{epsfig}
\usepackage{color}
\usepackage{tikz}
\usepackage[all,arc,dvips,arrow,curve,cmactex]{xy}

\usepackage{latexsym}
\usepackage{mathtools}
\DeclarePairedDelimiter\floor{\lfloor}{\rfloor}

\newtheorem{lemma}{Lemma}

\newtheorem{conj}{Conjecture}
\newtheorem{definition}{Definition}
\newtheorem{theorem}{Theorem}
\newtheorem{proposition}{Proposition}

\newcommand{\bea}{\begin{eqnarray}}
\newcommand{\eea}{\end{eqnarray}}
\newcommand{\beq}{\begin{equation}}
\newcommand{\eeq}{\end{equation}}
\newcommand{\enn}{\nonumber \end{equation}}

\title[$(f_0, f_1)$ pairs for convex polytopes]{A complete characterization of $(f_0, f_1)$-pairs of 6-polytopes}

\author{Karim Adiprasito}
\address[K.A.]{Einstein Institute of Mathematics, The Hebrew University of Jerusalem, Giv'at Ram, Jerusalem, 91904, Israel}
\email{karim.adiprasito@mail.huji.ac.il }

\author{R\'emi Cocou Avohou}
\address[R.C.A.]{
Max Planck Institut f\"ur Mathematik, Vivatsgasse 7, 53111 Bonn, Germany, \& ICMPA-UNESCO Chair, 072BP50, Cotonou,
Rep. of Benin
}
\email{avohou.r.cocou@mpim-bonn.mpg.de}

\begin{document}

\maketitle

\begin{abstract}
We completely characterize the first two entries, namely the $(f_0, f_1)$-vector pairs, for $6$-dimension polytopes. We also find the characterization for $7$-dimension polytopes with excess degree greater than $11$ and, we conjecture bounds fulfilled by $(f_0, f_1)$-vector pairs for any $d$-polytope having an excess degree greater than $3d - 10$.
%
%

\end{abstract}

\tableofcontents

\section{Introduction}
The study of $f$-vectors for convex polytopes have had huge successes in the last years. For a convex $d$-dimensional or $d$-polytope $P$, we write $f_i=f_i(P)$
for the number of $i$-dimensional faces (or $i$-faces) of $P$. The $f$-vector of $P$ is the vector $\big(f_0, f_1, \cdots, f_{d-1}\big)$ where faces of $P$ of dimension $0$, $1$, $2$, $d-2$ and $d-1$ are called vertices, edges, subfacets (or ridges), and facets of $P$, respectively. The number $2f_1-df_0:=\epsilon(P)$ is called the excess degree or excess of $P$.

In 1906 the $f$-vectors set of 3-polytopes was characterized by Steinitz \cite{Steinitz} and in 1980  Billera \& Lee and Stanley \cite{BilleraLee, Stanley} have proved the characterization of
the $f$-vectors of simplicial and of simple polytopes conjectured by McMullen \cite{PMcMullen} in 1971 through the famous ``$g$-theorem''. However little result seems to be known for general polytopes and the characterization of $f$-vectors of $d$-polytopes ($d\geq 4$) is still a big open problem in convex geometry. 

In most of the papers going in this direction, the authors have reduced the number of components of these vectors. For example Gr{\"u}mbaum, Barnette and Barnette-Reay \cite{BGrunbaum2, Barnette1, BarnetteReay} have characterized  for any $0\leq i<j\leq 3$ the following sets:
$$\Big\{\big(f_i, f_j\big): P \mbox{ is a $4$-polytope }\Big\}.$$

More recently, most of the results determine the possible number of edges or $1$-dimensional faces of a polytope, given the number of vertices.
Let $$\mathcal{E}^d=\left\{(f_0, f_1): \mbox{ there is a $d$-polytope with $f_0$ vertices and $f_1$ edges}\right\}.$$
Steinitz found the characterization for $d=3$, $\mathcal{E}^3=\left\{(f_0, f_1): \frac{3}{2}f_0\leq f_1\leq 3f_0-6\right\}.$ For $d=4, 5$ the results is given by the set $\left\{(f_0, f_1): \frac{d}{2}f_0\leq f_1\leq \binom{f_0}{2}\right\}$ from which some exceptions have been removed. Gr{\"u}mbaum \cite{BGrunbaum} proved the case $d=4$ by removing four exceptions: $(6, 12)$, $(7, 14)$, $(8, 17)$ and $(10, 20)$. The case $d=5$ becomes more complicated and has been proved in two different ways by G. Pineda-Villavicencio, J. Ugon and D. Yost, \cite{Guillermo2} and more recently by T. Kusunoki and S. Murai \cite{Kusunoki}. For this case, exceptions are infinitely many:
$$\mathcal{E}^5=\left\{(f_0, f_1): \frac{5}{2}f_0\leq f_1\leq \binom{f_0}{2}\right\}\setminus\left(\left\{\Bigg(f_0, \floor*{\frac{5}{2}f_0+1}\Bigg): f_0\geq 7\right\}\cup\Bigg\{(8, 20), (9, 25), (13, 35)\Bigg\}\right),$$
where $\floor*{r}$ denotes the integer part of a rational number $r$.

This paper addresses the study of $f$-vectors of convex polytopes and then thoroughly characterizes the set $\mathcal{E}^6$ i.e. the characterization for the first two entries of $6$ dimensional polytopes. The $(f_0, f_1)$-pairs characterizing the convex $7$-polytopes have been found under the condition that  $2f_1-7f_0>11$.
 Finally we conjecture bounds fulfilled by $(f_0, f_1)$ for arbitrary dimension $d> 4$.
 
 The  first part of this paper aims to find $\mathcal{E}^6$ by combining the ideas developed in \cite{Kusunoki, Guillermo2}. The main ingredients introduced in these two papers are reviewed in section \ref{sect:two}. Section \ref{sect:three} focuses on the proof of our first main result, i.e. Theorem \ref{theo:dimsix}, fully characterizes the first two entries of $6$-dimensional polytopes. Section \ref{sect:four} gives a complete characterization of the $(f_0, f_1)$-pairs of $7$-polytopes with excess greater than $11$  and addresses our second main result in Theorem \ref{theo:dimsept}. An extension of this result in arbitrary dimensions is conjectured in the same section \ref{sect:four} (see Conjecture \ref{conj:dimarb}).


\section{Preliminaries}\label{sect:two}
This subsection recalls some recent results on face numbers of convex polytopes. Let a polytope $Q$ be the pyramid over a polytope $P$. Then we have
$$f_0(Q)=f_0(P)+1, \quad f_1(Q)=f_0(P)+f_1(P).$$
If the polytope $P$ is of dimension $d$ and has a simple vertex $v$ and $Q$ is obtained from $P$  by truncating the vertex $v$ then $$f_0(Q)=f_0(P)+d-1\mbox{ and }  f_1(Q)=f_0(P)+\binom{d}{2}.$$ 
These relations are the most important in the proof of the following theorem.
\begin{theorem}\label{theo:fivech}\cite{Kusunoki} The set of $(f_0, f_1)$-vector pairs for $5$-polytopes is given by
$$\mathcal{E}^5=\left\{(f_0, f_1): \frac{5}{2}f_0\leq f_1\leq \binom{f_0}{2}\right\}\setminus\left(\left\{\Bigg(f_0, \floor*{\frac{5}{2}f_0+1}\Bigg): f_0\geq 7\right\}\cup\Bigg\{(8, 20), (9, 25), (13, 35)\Bigg\}\right).$$
\end{theorem}

Another quantity called the excess degree of a polytope is introduced in \cite{Guillermo2} and proves useful to simplify the above theorem.
\begin{definition}[Excess degree]
The excess degree or excess $\Sigma(P)$  of a $d$-polytope $P$ is defined as the sum of the excess
degrees of its vertices and given by $\epsilon(P)=2f_1-df_0$.
\end{definition}

We now review the polytopes with small excess degrees starting from the simple polytopes to the $d$-polytopes with excess $d-2$ and $d-1$.
\begin{theorem}\label{theo:exess}\cite{Guillermo2}
Let $P$ be a $d$-polytope. Then the smallest values in $\Sigma(d)$ are $0$ and $d-2$.
\end{theorem}

\begin{theorem}\label{theo:combeq}\cite{Guillermo}
Up to combinatorial equivalence,

\begin{itemize}
\item[(i)] the simplex $\Delta_{0, d}$ is the only simple $d$-polytope with strictly less than $2d$ vertices;
\item[(ii)] the simplicial prism $\Delta_{1, d-1}$ is the only simple $d$-polytope with between $2d$ and $3d-4$ vertices;
\item[(iii)] $\Delta_{2, d-2}$ is the only simple $d$-polytope with $3d-3$ vertices;
\item[(iv)] the only simple $d$-polytope with $3d-2$ vertices is the $6$-dimensional polytope $\Delta_{3, 3}$;
\item[(v)] the only simple $d$-polytope with $3d-1$ vertices are the polytope $J_d$, the $3-$dimensional cube $\Delta_{1, 1, 1}$ and the $7-$dimensional polytope $\Delta_{3, 4}$;
\item[(v)] there is a simple polytope with $3d$ vertices if, $d=4, 8$; the only possible examples are $\Delta_{1, 1, 2}$, $\Gamma_{2, 2}$ and  $\Delta_{3, 5}$.
\end{itemize}
\end{theorem}

\begin{theorem}\label{theo:exessdmoindeux}\cite{Guillermo}
Any $d$-polytope $P$ with excess exactly $d-2$ either 

\begin{itemize}
\item[(i)] has a unique nonsimple vertex, which is the intersection  of two facets, or
\item[(ii)] has $d-2$ vertices of excess degree one, which form a $d-3$-simplex which is the intersection of two facets.
\end{itemize}
In either case, the two intersecting facets are both simple polytopes, $P$ is a Shepard polytope and has another facet with excess $d-3$.
\end{theorem}

The following statement characterizes the $d$-polytopes excess degree $d-1$.

\begin{theorem}\label{theo:exessdmoinun}\cite{Guillermo}
Let $P$ be a $d$-polytope $P$ with excess exactly $d-1$, where $d>3$. Then $d=5$ and either 

\begin{itemize}
\item[(i)] there is a single vertex with excess four, which  is the intersection of two facets, and whose vertex figure is $\Delta_{2, 2}$; or
\item[(ii)] there are two vertices with excess two, the edge joining them is the intersection of two facets and its underfacet is either $\Delta_{1, 1, 2}$ or $\Gamma_{2, 2}$; or
\item[(iii)] there are four vertices each with excess one, which form a quadrilateral $2$-face which the intersection os two facets, and whose undefacet is the tesseract $\Delta_{1, 1, 1, 1}$.
\end{itemize}
In all cases $P$ is a Shephard polytope.
\end{theorem}

We set for all $d$-dimensional polytopes $\phi(v, d)=\frac{1}{2}dv+\frac{1}{2}(v-d-1)(2d-v)$.
\begin{theorem}\label{theo:cond}\cite{Guillermo2}
Let $P$ be a $d$-polytope.
\begin{enumerate}
\item If $f_0(P)\leq 2d$, then $f_1(P)\geq \phi(f_0(P), d)$.
\item If $d\geq 4$, then $(f_0(P), f_1(P))\neq (d+4, \phi(d+4, d)+1)$.
\end{enumerate}
\end{theorem}

\begin{theorem}\label{theo:nottriplex}\cite{Guillermo2}
If $4\leq k \leq d$, then a $d$-polytope with $v = d + k$ vertices which is not a triplex must have at least $\phi( v, d ) + k - 3$ edges. In other words, its excess degree is at least $( k - 1 )( d - k ) + 2 ( k - 3 )$.
\end{theorem}

We now recall the $g$-conjecture for simplicial complex  \cite{McMullen}. For a simplicial complex $\Delta$ of dimension $d$, its $f$-vector is $(f_0(\Delta), \cdots, f_d(\Delta))$ where $f_i(\Delta)$ is the number of $i$-dimensional faces of $\Delta$. The $h$-vector is $(h_0(\Delta), \cdots, h_d(\Delta))$ where 
$$h_k=\sum_{i=0}^k(-1)^{k-i} \binom{d+1-i}{d+1-k}f_{i-1};$$ $\forall$ $k=0, \cdots, d+1$.
\begin{theorem}\label{gtheor} \cite{McMullen}
Let $\Delta$ be a simplicial $d$-sphere. A sequence of integers $(h_0, \cdots, h_d)$ is the $h$-vector of $\Delta$ if and only if the following two conditions hold:
\begin{enumerate}
\item[i)] $h_i=h_{d-i}$ for all $0\leq i\leq d$.
\item[ii)] $\Big(h_0, h_i-h_2,  \cdots, h_{\floor*{\frac{d+1}{2}}}-h_{\floor*{\frac{d+1}{2}}-1}\Big)$ is an $M$-vector.
\end{enumerate}
\end{theorem}

The following theorem in \cite{G.kalai2} was conjectured in \cite{G.kalai1}, and extends the lower bound theorem to arbitrary $d$-polytopes. Let $C$ be a polyhedral complex and $f_2^k(C)$ the number of $2$-faces of $C$ which are $k$-gons. 

\begin{theorem}\label{theokalai}
If P is a $d$-polytope with $n$ vertices then 
$$f_1(P)+\sum_{k\geq 3}(k-3)f_2^k(P)\geq d\times n-\binom{d+1}{2}.$$
\end{theorem}

The graph $G(P)$ of a polytope $P$ is the abstract graph with vertex set $V(P)$ and edge set $E(P)$ denoting respectively the set of vertices and edges of $P$. This section is ended with the Balinski's theorem which claims that the graph of a $d$-polytope is $d$-connected.

\begin{theorem}\cite{Balinski}
Given a $d$-dimensional convex polytope $P\subset \mathbb{R}^N$, the graph $G(P)\setminus V$is connected for any subset $V$ of the vertex set of $P$ which is contained in some $(d-2)$-dimensional affine subspace of $\mathbb{R}^N$. In particular, $G(P)$ is $d$-connected.
\end{theorem}

\section{Complete characterization of all possible $(f_0, f_1)$-pairs of $6$-polytopes}\label{sect:three}

If $P$ is a $6$-polytope having a simple vertex $v$ and $Q$ the $6$-polytope obtained from $P$  by truncating the vertex $v$ then
$$f_0(Q)=f_0(P)+5\mbox{ and }  f_1(Q)=f_0(P)+15.$$ We can prove that if for a $6$-polytope $P$ we have $f_1(P)\leq \frac{7}{2}f_0$ then $P$ has at least one simple vertex.

The main result in this section is the following
\begin{theorem} \label{theo:dimsix}The set of $(f_0, f_1)$-vectors pair for $6$-polytopes is given by
\begin{eqnarray}
\mathcal{E}^6=\left\{(f_0, f_1): 3f_0\leq f_1\leq \binom{f_0}{2}\right\}\setminus \Bigg(\left\{\Big(f_0, 3f_0+1\Big): f_0\geq 7\right\}\cup&&\cr \Big\{(8, 24); (9, 27); (9, 29); (10, 30); (10, 32); (10, 34); (11, 33);  (12, 38); (12, 39); (13, 39); (14, 42);&&\cr  (14, 44); (15, 47); (18, 54); (19, 57);(17, 53); (20, 62)\Big\}\Bigg).
\end{eqnarray}
\end{theorem}
To prove this result we must establish the following propositions and lemmas.

\begin{proposition} \label{Prop:parity}
If $P$ is a $d$-polytope with $d> 4$, then
\begin{eqnarray}
f_1(P)\neq \floor*{\frac{d}{2}f_0(P)+1}.
\end{eqnarray}
\end{proposition}

\proof
Assume that $f_1(P)= \floor*{\frac{d}{2}f_0(P)+1}$. If $f_0(P)$ is even then  $2f_1(P)-df_0(P)=2$ and $0<2<d-2$ which is impossible since  Theorem \ref{theo:exess} states that $\Sigma(d)$ can not take any value between $0$ and $d-2$. If $f_0(P)$ is odd then $2f_1(P)-df_0(P)=1$ and $0<1<d-2$ which is also impossible.
\qed

In particular for $d=6$, $(f_0, 3f_0+1)\notin \mathcal{E}^6$ and then $(f_0+5, 3f_0+1+15)=(f_0+5, 3(f_0+5)+1)\notin \mathcal{E}^6$.

\begin{lemma}\label{lemma:neufdix}
The following relations hold $(8, 24); (9, 27); (9, 29); (10, 30); (10, 32); (10, 34); (11, 33)\notin \mathcal{E}^6$ and $\big(f_0+1, \floor*{\frac{7}{2}f_0+1}\big)\notin \mathcal{E}^6$ for $f_0=7, 8, 9$.
\end{lemma}
\proof
The fact that $(10, 34)\notin \mathcal{E}^6$ is given by Theorem \ref{theo:cond}
(2) and all the remaining are given by Theorem \ref{theo:cond} (1).
\qed

The following lemma is obtained from pyramids over $5$-polytopes in Theorem \ref{theo:fivech}.
\begin{lemma}\label{lemma:fpyra}
This inclusion is true:
\bea
\left\{(f_0, f_1): \frac{7}{2}f_0-\frac{7}{2}\leq f_1\leq \binom{f_0}{2}\right\}\setminus\cr\left(\left\{\Bigg(f_0+1, \floor*{\frac{7}{2}f_0+1}\Bigg): f_0\geq 7\right\}\cup\Bigg\{(9, 28), (10, 34), (14, 48)\Bigg\}\right)\subset \mathcal{E}^6.
\eea

\end{lemma}
Let us discuss the set $\left\{\Bigg(f_0+1, \floor*{\frac{7}{2}f_0+1}\Bigg): f_0\geq 7\right\}$. Lemma \ref{lemma:neufdix} already addressed the cases $f_0=7, 8, 9$.

\begin{lemma}\label{lemma:trois}
\begin{itemize}
\item There is no $6$-polytope with $11$ vertices and $36$ edges and no $6$-polytope with $12$ vertices and $38$ edges. 
\item The following exist $(13, 43), (14, 48)\in \mathcal{E}^6$.
\end{itemize}
\end{lemma}
\proof
The $6$-polytopes with $11=6+5$ vertices and $36$ edges or the $6$-polytopes with $12=6+6$ vertices and $38$ edges are not triplexes. Their excess degrees  are respectively $6$ and $4$ but from Theorem \ref{theo:nottriplex}, their excess degrees should be at least $8$ and $6$ respectively. Then $(11, 36), (12, 38) \notin \mathcal{E}^6$.

Consider the duals of the cyclic $6$-polytopes $C$ and $C'$ with $7$ and $8$ vertices respectively. Performing the connected sum between these two polytopes we obtain a polytope with $13 = 8 + 7 - 2$ facets and  $43= 28 + 21 - 6$ codimension 2-faces giving $(13, 43)\in \mathcal{E}^6$.
%

Consider the cyclic $6$-polytope $C$ with $7$ vertices and let us add a pyramid over one of its facets. We obtain a polytope $P$ whose $(f_0, f_1)$ pair is $(8, 27)$. The truncation of a simple vertex of $P$ give a polytope $P'$ with $(f_0, f_1) = (13, 42)$. As the truncation of a simplex vertex create a simplex facet then $P'$ contains a simplex facet  $F$. Adding a pyramid over $F$ gives a polytope whose $(f_0, f_1)$ pair is $(14, 48)$.
\qed

The case $(12, 39)\notin \mathcal{E}^6$ will not need any proof if the conjecture in \cite{Guillermo2} claiming that for $d\geq 6$, there are no $d$-polytopes with $2d$ vertices and $d^2+d-3$ edges was true but we have the following.

\begin{lemma}\label{lemma:fseptc}
There is no $6$-polytope with $12$ vertices and $39$ edges. 
\end{lemma}
\proof
We assume that such polytope exists and is denoted by $P$. Let $F$ be a facet of $P$.

The $5$-polytope $F$ has at most $11$ vertices and the particular case of $11$ vertices is not possible. In fact $F$ would have at least 29 edges as there is no $5$-polytope with $11$ vertices and $\floor*{\frac{5}{2}\times 11+1}=28$ edges. This will make $P$ to have at least $40$ edges. 

Suppose that $F$ has $10$ vertices. Its number of edges can not exceed $27$ since there should be at least $12$ edges outside $F$. Therefore the facet $F$ has $25$ or  $27$ edges because there is no $5$-polytope with $11$ vertices and $26$ edges. If the number of edges in $F$ is $25$ then $F$ is the simplicial prism $\Delta_{1, 4}$. Let $F'$ be the other facet is $P$ intercepting $F$ at the ridge $\Delta_{0, 4}$. The $5$-polytope $F'$ has at most $6$ vertices and therefore is the $5$-simplex. The set $F\cup F'$ has $11$ vertices and $30$ edges causing $P$ to have $42$ edges, which is a contradiction. In the case where the $5$-polytope $F$ has $27$ edges, then its excess degree is $5-1=4$. From Theorem \ref{theo:exessdmoinun} we much evaluate three consecutive cases: assume that there is a single vertex with excess four, which is the intersection of three facets. As the minimum number of vertices of each of these three facets $F_i$; $i=1, 2, 3$ is $5$ then the total number of vertices in $\cup_{i}F_i$ exceeds $10$, which is a contradiction. Another possibility is that there are two vertices with excess two, the edge joining them is the intersection of two facets. Let us denote them by $F_1$ and $F_2$. The minimum number of vertices in $F_i$, $i=1, 2$ is $5$ and then $F_1\cup F_2$ has at least $8$ vertices and $19$ edges. But we have at least $9$ edges going outside $F_1\cup F_2$ and touching the two remaining vertices of $F$. This is impossible as it will cause $F$ to have at least $28$ edges. We now consider the last case of Theorem \ref{theo:exessdmoinun}  saying that there are four vertices each with excess one, which form a quadrilateral $2$-face which is the intersection of two facets. Having a quadrilateral $2$-face, these two facets are not simplices and each of them have at least $6$ vertices and $13$ edges since they is no $4$-polytope with $6$ vertices and $12$ edges. Therefore the union of these two facets will have at least $8$ vertices and $22$ edges. Adding the $22$ edges to the edges joining the two remaining vertices of $F$ leads to a contradiction.

Let us now study the case where the facet $F$ has $9$ vertices. The minimum number of edges in $F$ is $24$ and should be exactly $24$ since we already have at least $15$ edges outside $F$. Its excess degree $5-2=3$ implies from Theorem \ref{theo:exessdmoindeux} that: $F$ has a unique non-simple vertex which is the intersection of two simple facets $R_1$, $R_2$ (ridges of $P$) or $F$ has two simple facets which intercept at a $3$-simplex. Each of the $4$-polytopes $R_1$ and $R_2$ has at least $5$ vertices and $R_1\cup R_2$ has at least $9$ vertices which is impossible. The second possibility is that $F$ has two simple facets denoted again $R_1$ and $R_2$ which intercept at a $3$-simplex. In this case $R_1\cup R_2$ has at least $6$ vertices and $14$ edges. This is impossible since the number of edges touching at least one vertex in $F\setminus (R_1\cup R_2)$ is at least $12$ causing $F$ to have at least $26$ edges. 
 
In the case where the facet $F$ has $8$ vertices there is also a contradiction. Indeed the number of edges should be at least $22$. Furthermore the number of edges outside $F$ is at least $18$ and $P$ should have at least $40$ edges.

If the number of vertices in $F$ is $7$, then its number of edges is necessarly $19$ since we have at least $20$ edges in $P$ which are outside $F$. Its excess degree is $5-2=3$ and proceeding in the same way as above leads to a contradiction. 

Finally let us assume that every facet $F$ of $P$ has $6$ vertices. This implies that$F$ is a simplex and hence a $2$-simplicial polytope. The polytope $P$ is therefore $2$-simplicial and Theorem \ref{theokalai} claims that it would have at least $51$ edges, giving the last contradiction.
\qed

\begin{lemma}\label{lemma:fsept}
For an odd integer $f_0\geq 12$ we have $\big(f_0+1, \floor*{\frac{7}{2}f_0+1}\big)\in \mathcal{E}^6$. Furthermore if $\big(f_0+1, \floor*{\frac{7}{2}f_0+1}\big)\in \mathcal{E}^6$, then $\big(f_0+7, \floor*{\frac{7}{2}(f_0+6)+1}\big)\in \mathcal{E}^6$.
\end{lemma}
\proof
Suppose that $f_0$ is odd. Since $f_0\geq 12$ then $f_0-4\geq 8$ and from Lemma \ref{lemma:fpyra}, $(f_0-4, \floor*{\frac{7}{2}(f_0-4)})\in \mathcal{E}^6$. Also $\floor*{\frac{7}{2}(f_0-4)}<\frac{7}{2}(f_0-4)$ as $f_0-4$ is odd and then $\big(f_0+1, \floor*{\frac{7}{2}f_0+1}\big)\in \mathcal{E}^6$  by truncation of simple vertex.

Let $P$ be a $6$-polytope with $(f_0, f_1)$-pairs equal to $\big(f_0(P)+1, \floor*{\frac{7}{2}f_0(P)+1}\big)\in \mathcal{E}^6$ then after the truncation of a simple vertex of $P$ and a pyramid over a simplex facet of the resulting polytope we obtain a $6$-polytope $Q$ with $f_0(Q)=f_0(P)+7$ and $f_1(Q)=\floor*{\frac{7}{2}(f_0(P)+6)+1}$.
\qed

\begin{lemma} For any integer $f_0$ satisfying $f_0\geq 12$, 
 $\big(f_0+1, \floor*{\frac{7}{2}f_0+1}\big)\in \mathcal{E}^6$.
\end{lemma}
\proof
Assume that $f_0\geq 12$. From the second result of Lemma \ref{lemma:fsept} it is enough to check the result for $f_0=12, 13, 14, 15, 16, 17$. The case $f_0=12$ comes from Lemma \ref{lemma:fsept}. Applying the first result in Lemma \ref{lemma:fsept} we can conclude the result for $f_0=13, 15, 17$. Let us now focus on the two remaining cases $f_0=14, 16$ which are $(15, 50)$ and $(17, 57)$.

Consider the $6$-polytope $P$ with $f_0(P)=10$ and $f_1(P)=35$ obtained from a pyramid over a $5$-polytope $Q$. If we assume that $P$ has no simple vertex then each of its vertices has degree $7$ since $\sum_{v\in P} deg(v)=70$ and this is impossible since taking a pyramid over $Q$ implies that $P$ has a vertex of degree $9$. Then $P$ has a simple vertex which truncation gives a polytope $P'$ with $f_0(P')=15$ and $f_1(P)=50$. Hence $(15, 50)\in \mathcal{E}^6$.

Let $R$ be a $6$-polytope with $f_0(R)=12$ and $f_1(P)=42$ obtained from a pyramid over a $5$-polytope. The same procedure as above gives $(17, 57)\in \mathcal{E}^6$.
\qed


\begin{lemma}\label{lemma:extracases}
\begin{itemize}
\item The following polytopes pairs do not exist: $$(13, 39); (14, 42); (14, 44); (15, 47); (18, 54), (19, 57)\notin \mathcal{E}^6.$$
\item The following pairs are possible: $$(15, 45); (15, 49); (16, 48); (17, 54); (19, 59); (23, 69); (24, 72); (27, 83); (35, 107)\in \mathcal{E}^6.$$
\end{itemize}
\end{lemma}
\proof
Assume that there is a $6$-polytope $P$ with $13=2\times 6+1$ vertices and $39$ edges. The polytope $P$ is not a pentasm as the number of edges is smaller then $41$ and this contradict the fact that for $d\geq 5$ there is no polytope with $2d+1$ vertices and fewer than $d^2+d-1$ edges. Hence $(13, 39)\notin \mathcal{E}^6$.
Knowing that $\min E(14, 6)=45$ then $(14, 42); (14, 44)\notin \mathcal{E}^6$. Let us assume that there is a $6$-polytope $Q$ with $f_0(Q)=15$ and $f_1(Q)=47$. The excess of $Q$ is $4=6-2$. From Theorem \ref{theo:exessdmoindeux} the polytope $Q$ has $4$ vertices of excess degree one, which form a $3$-simplex or $Q$ has a unique non-simple vertex of degree 10. Suppose that $Q$ has 4 vertices forming a $3$-simplex. The same theorem says that this $3$-simplex is in the intersection of two simple facets. Let $F$, $G$ be these two facets. If they both have less than 10 vertices then $F=\Delta_{0, 5}=G$. The number of edges in $F\cup G$ is $30-6=24$ and the number of edges incident to the set $S$ of the $7$ vertices outside  $F\cup G$ must be $23$. If $x$ is the number of edges between the 7 vertices and $y$ is the number of edges between $F\cup G$ and $S$ then $x+y=23$ and $2x+y\geq 7\times 6$. This implies that $y\leq 4$ which means that removing at most four vertices in $S$ disconnects the polytope graph. This is a contradiction to the theorem of Balinski's. If $F$ has $10$ vertices and $G$ has 6 then $G=\Delta_{0, 5}$, $F=\Delta_{1, 4}$ and $F\cup G$ has 12 vertices, 34 edges and we need at least $6\times 3-3=15$ edges joining at least one of the three vertices outside $F\cup G$. This gives $49>47$ which is a contradiction. If $F$  has $12$ vertices and $G$ still has $6$ vertices then $G=\Delta_{0, 5}$ and $F=\Delta_{2, 3}$.
The intersection of $F$ and $G$ is the $3$-simplex whose vertices have degree 7 in 
 $F\cup G$ and then adding a simple vertex can not be connected to any of them. The removal of the $4$ vertices and the only vertex outside $F\cup G$ will disconnect the graph of $Q$ which is an absurdity from Balinski's theorem.  From Theorem \ref{theo:combeq} (iv) there is no simple $5$-polytope with $13$ vertices. If $G=\Delta_{0, 5}$ then $F$ can not have more than $13$ vertices. Otherwise $F\cup G$ will contain at least $16$ vertices. If $F$ and $G$ have at least $10$ vertices then  $F\cup G$ will contain at least $16$ vertices which is also impossible.

If $Q$ has a unique non-simple vertex $v$ then $deg(v)=10$ and $v$ is the intersection of two facets $F_1$ and $F_2$ which are two simple $5$-polytope from Theorem \ref{theo:exessdmoindeux}(i). If both of them have lest than $10=2\times 5$ vertices then from Theorem \ref{theo:combeq}(i) $F_1=\Delta_{0, 5}=F_2$. The number of edges in $F_1\cup F_2$ is 30 and the number of edges outside  $F_1\cup F_2$ should be at least $6\times 4-6=18$ and the sum is 48 which is impossible. Suppose that $F_2=\Delta_{0, 5}$ and $F_1$ has 10 vertices then $F_1=\Delta_{1, 4}$. The number of vertices in $F_1\cup F_2$ is $15$ and in oder to find $Q$ we need to add $7$ edges and no vertex. The simplex $\Delta_{0, 5}$ intercept the simplicial prism $\Delta_{1, 4}$ at $v$ and the remaining vertices in $\Delta_{0, 5}$ have degree $5$. If we want to add five edges to $F_1\cup F_2$ to increase the degrees of these vertices to $6$, some of these edges will be incident to vertices of the second simplex $\Delta_{0, 4}$ of $\Delta_{1, 4}$ which do not touch $\Delta_{0, 5}$.The removing of $v$ and these vertices (whose number is at most $4$) will disconnect the polytope graph. This is a contracdition.

From Theorem \ref{theo:combeq} (vi) there is a simple polytope with $3d$ vertices if and only if $d=4, 8$. Then $(18, 54)\notin \mathcal{E}^6$.
Let us now concentrate on the case $(19, 57)$ which is a simple polytope if it exists. Assume that there is a simple $6$-polytope  $P$ with $19$ vertices and $57$ edges. The dual $P^\star$ of $P$ is a simplicial polytope with $f$-vector sequence $(f_0, f_1, f_2, f_3, f_4, f_5)$ where $f_4=57$ and $f_5=19$. For all $d$-dimensional simplicial polytope the following inequality holds: $f_{d-1}\geq (d-1)f_0-(d+1)(d-2)$ \cite{Barnette}. Then $f_5\geq 5f_0-28$ implies that $f_0=8$ or $f_0=9$. The $g$-theorem (Theorem \ref{gtheor}) for simplicial polytopes says that the sequence of integers $(h_0, \cdots, h_7)$ is the $h$-vector of $P^\star$ with $$h_k=\sum_{i=0}^k(-1)^{k-i} \binom{d+1-i}{d+1-k}f_{i-1};$$ $\forall$ $k=0, \cdots, 7$. From the same theorem we also have $h_i=h_{7-i}$ $\forall$ $i=0, \cdots, 7$ and now compute the numbers $h_i's$ and obtain:
\bea
h_1&=&-7+f_0,\cr
h_2&=&21-6f_0+f_1,\cr
h_3&=&-35+15f_0-5f_1+f_2,\cr
h_4&=&35-20f_0+10f_1-4f_2+f_3,\cr
h_5&=&-21+15f_0-10f_1+6f_2-3f_3+f_4,\cr
h_6&=&7-6f_0+5f_1-4f_2+3f_3-2f_4+f_5.
\eea 
From $h_1=h_6$ and $h_2=h_5$ we get $f_2=\frac{1}{2}(28-14f_0+6f_1+f_4-f_5)$ and the system of equations $h_3=h_4$; $h_2=h_5$ also gives $f_2=\frac{1}{9}(168-84f_0+34f_1+f_4)$. Equaling these two expressions of $f_2$ we get $f_1=\frac{1}{14}(-84+42f_0+7f_4-9f_5)$ which is not an integer for $f_0=8, 9$. In conclusion $(19, 57)\notin \mathcal{E}^6$.

The $(f_0, f_1)$-pairs of the polytopes $\Delta_{2, 4}$ and $\Delta_{3, 3}$ are respectively $(15, 45)$ and $(16, 48)$. Taking a pyramid over the $5$-polytopes with $(f_0, f_1)$-pairs $(14, 35)\in \mathcal{E}^5$, we obtain a $6$-polytope which two first entries are $(15, 49)\in \mathcal{E}^6$ and truncating a simple edge from $6$-polytopes having $(9, 30)\in \mathcal{E}^6$, $(11, 35)\in \mathcal{E}^6$  and $(15, 45)\in \mathcal{E}^6$ as $(f_0, f_1)$-pairs give respectively the $6$-polytopes whose two first entries are respectively $(17, 54)\in \mathcal{E}^6$, $(19, 59)\in \mathcal{E}^6$ and $(23, 69)\in \mathcal{E}^6$. Similarly the truncation of a simple edge from $\Delta_{3, 3}$ and the $6$-polytope with $(f_0, f_1)$-pair equal to $(19, 59)\in \mathcal{E}^6$ give respectively the $6$-polytopes with $(24, 72)\in \mathcal{E}^6$ and  $(27, 83)\in \mathcal{E}^6$ as $(f_0, f_1)$-pairs. The same operation on  $(27, 83)$ gives $(35, 107)\in \mathcal{E}^6$.
\qed

Let us now concentrate on the set $f_1\in \{3f_0\}\cup\big]3f_0+1, \frac{7}{2}f_0-\frac{7}{2}\big[$ for $f_0\geq 7$. We set:
\begin{eqnarray}
X'=\Big\{(8, 24); (9, 27); (9, 29); (10, 30); (10, 32); (10, 34); (11, 33); (11, 36); (12, 38); (12, 39); (13, 39);&&\cr (14, 42);  (14, 44); (15, 47); (18, 54); (17, 53); (19, 57); (20, 62)\Big\}
\end{eqnarray}

\begin{lemma}
For $f_0\geq 7$; if $(f_0, f_1)\notin X'$ and $f_1\in \{3f_0\}\cup\big]3f_0+1, \frac{7}{2}f_0-\frac{7}{2}\big[$ then $(f_0, f_1)\in \mathcal{E}^6$.
\end{lemma}
\proof
Assume that $7\leq f_0< 13$ and $f_1\in \{3f_0\}\cup\big]3f_0+1, \frac{7}{2}f_0-\frac{7}{2}\big[$. In this case $f_1=3f_0$. For $f_0=7$ we have $f_1=21$ which is known as the cyclic $6$-polytope with 7 vertices. If $f_0=8, 9, 10, 11,  12$ then we have respectively $f_1=24$; $f_1\in\{27, 28\}$; $f_1\in\{30, 31\}$; $f_1\in\{33, 34\}$ and $f_1\in\{36, 37, 38\}$. The fact that the simplicial prism $\Delta_{1, 5}$ is a $6$-polytope having $12$ vertices and $36$ edges gives $(12, 36)\in \mathcal{E}^6$ . The remaining does not exist and the proofs come from Proposition \ref{Prop:parity} and Lemma \ref{lemma:neufdix}. 

Let us now discuss the cases $f_0\geq 13$. If  $f_1=3f_0$ then let $k$ be the biggest integer such that $k< \frac{1}{4}\Big(7f_0-7-2f_1\Big)$. Setting $f_0'=f_0-5k$ and $f_1'=f_1-15k$; we have $f_1'=3f_0'$. If $(f_0', f_1')\in X'$ then $(f_0', f_1')=(8, 24); (9, 27); (10, 34); (11, 33); (13, 39); (14, 42); (18, 54)$. We concentrate first on the following cases: $(f_0', f_1')=(8, 24); (13, 39); (18, 54)$. It was proved in Lemma \ref{lemma:extracases} that $(23, 65)\in \mathcal{E}^6$ and then  $(23+5k, 69+15k)\in \mathcal{E}^6$ by $k$ repeated truncation of simple vertices and also the unfeasibility of  $(8, 24); (13, 39)$ and $(18, 54)$ have been proved. In the same way we have $(f_0', f_1')=(9, 27); (14, 42); (19, 57)\notin \mathcal{E}^6$ and from Lemma \ref{lemma:extracases}  $(24, 72)\in \mathcal{E}^6$ and then $(24+5k, 72+15k)\in \mathcal{E}^6$ for suitable $k$. Also $(f_0', f_1')=(10, 30); (11, 33)\notin \mathcal{E}^6$ and we have $(15, 45); (16, 48)\in \mathcal{E}^6$ which imply that $(15+5k, 45+15k), (16+5k, 48+15k)\in \mathcal{E}^6$.

For $f_1\in \big]3f_0+1, \frac{7}{2}f_0-\frac{7}{2}\big[$ let $k$ be the smallest integer such that $k\geq \frac{1}{4}\Big(7f_0-7-2f_1\Big)$. We set $f_0'=f_0-5k$; $f_1'=f_1-15k$. The following inequality holds $\frac{7}{2}f_0'-\frac{7}{2}\leq f_1'<\frac{7}{2}f_0'$. If $(f_0', f_1')\notin X'$ then $(f_0', f_1')\in \mathcal{E}^6$ and truncating a simple vertex $k$ times gives a polytope with $f_0$ vertices and $f_1$ edges. Otherwise $(f_0', f_1')=(9, 29); (10, 32); (10, 34); (12, 38); (12, 39); (14, 44); (15, 47); (17, 53); (20, 62)$; $(22, 68); (25, 77); (30, 92)$. It is known that $(9, 29); (14, 44)\notin \mathcal{E}^6$. The case $(19, 59)\in \mathcal{E}^6$ was solved in Lemma \ref{lemma:extracases} and all the cases $(19+5k, 59+15k)$ are obtained by successive truncation of simple vertices. We now suppose that $(f_0', f_1')=(10, 32); (15, 47); (20, 62); (25, 77); (30, 92)$. From Lemma \ref{lemma:extracases} we proved that $(35, 107)\in \mathcal{E}^6$ and then $(35+5k, 107+15k)\in \mathcal{E}^6$ by performing $k$ repeated truncation of simple vertices. We have also proved the case $(15, 47)\in \mathcal{E}^6$. It remains to evaluate the cases $(20, 62); (25, 77); (30, 92)$. For $(f_0', f_1')=(10, 34)$ or  $(f_0', f_1')=(12, 39)$; we have proved that $(15, 59); (17, 54)\in \mathcal{E}^6$ and with the same argument as above $(15+5k, 59+15k); (17+5k, 54+15k)\in \mathcal{E}^6$ for suitable $k$. Finally 
let us study the case $(f_0', f_1')=(12, 38); (17, 53); (22, 68)$. It was proved that $(27, 83)\in \mathcal{E}^6$ and then $ (27+5k, 83+15k)\in \mathcal{E}^6$. The case $(12, 38)$ is given by Lemma \ref{lemma:trois} and the remaining cases are $(17, 53); (20, 62); (22, 68); (25, 77); (30, 92)$. 
\qed

The following lemma solves the remaining cases. 
\begin{lemma}\label{lemma:moreextracases}
\begin{itemize}
\item The cases $(17, 53); (20, 62)$ are unfeasible.
\item This relation $(22, 68); (25, 77); (30, 92) \in \mathcal{E}^6$ holds.
\end{itemize} 
\end{lemma}

\proof
We start with the case $(17, 53)$. If there is a $6$-polytope $P$ with $17$ vertices and $53$ edges then its excess degree is $6-2=4$. Theorem \ref{theo:exessdmoindeux} claims that $P$ has a unique non-simple vertices which is the intersection of two simple facets or $P$ has four non-simple vertices forming a $3$-simplex and which is in the intersection of two simple facets. In case we have a unique non-simple vertex $v$ each of the two simple facets $F_1$, $F_2$ intersecting at $v$ has at least $6$ vertices and $10$ edges. There union will have at least $11$ vertices and $30$ edges. The same theorem states that $P$ has another facet $F_3$ of excess $6-3=3$. This facet intercepts $F_1$ at a ridge $R_1$ and $F_2$ at a ridge $R_2$. Therefore $F_3$ has at least  $9$ vertices which is the minimum number of vertices in $R_1\cup R_2$ and $24$ edges.  The minimum number of vertices in $F_1\cup F_2\cup F_3$ is $11$ and we have at least $30+24-20=34$ edges. This is impossible since the number of edges touching the remaining $6$ vertices of $P$ is at least $21$ causing $P$ to have at least $55$ edges. We now consider the second case of $4$ non-simple vertices. Let us denote by $F_1$ and $F_2$ the two facets intercepting at the $3$-simplex formed by the four vertices of excess one.  The same analyze as earlier helps us to find out that the minimum number of vertices in $F\cup F_2$ is $8$ and we have at least $24$ edges. Let $e_a$ be the number of edges between the $9$ remaining vertices and $e_b$ be the number of edges  joinning any of the $9$ vertices and a vertex in $F \cup F_2$. It comes out that $e_a+e_b=29$ and $2e_a+e_b\geq 6\times 9$. This implies that $e_a\geq 25$ and $e_b\leq 4$ which contradicts Balinski's theorem. In conclusion $(17, 53) \notin \mathcal{E}^6$.

Let us study the case $(20, 62)$ by assuming that the polytope $P$ exists. Its excess degree is $3=5-2$ and we can proceed as in the case $(17, 53)$. The case of a unique non-simple vertex leads to the facets $F_1$, $F_2$ and $F_3$ in which we have again $11$ as the minimum number of vertices in $F_1\cup F_2\cup F_3$ and at least $30+24-20=34$ edges. Denoting by $e_a$ be the number of edges between the $9$ remaining vertices and $e_b$ be the number of edges  joining any of the $9$ vertices and a vertex in $F_1\cup F_2\cup F_3$ we have the same relations $e_a+e_b=29$ and $2e_a+e_b\geq 6\times 9$. This is also run out by Balinski's theorem. We now turn to the case of two simple facets $F_1$ and $F_2$ intercepting at a $3$-simplex. As $F_1$ and $F_2$  are simple then the minimum number of vertices and edges in each of them is respectively $6$ and $15$. Therefore the number of vertices in $F\cup F_2$ is at least $8$ with a minimum of $24$ edges. Proceeding in the same way as above we obtain the relations  $e_a+e_b=38$ and $2e_a+e_b\geq 6\times 12$ since there are at least $12$ vertices outside $F_1\cup F_2$ and at least $38$ edges joining them. This also contradicted by Balinski's theorem and then $(20, 62) \notin \mathcal{E}^6$.

We now turn to the cases $(22, 68)$, $(25, 77)$ and $(30, 92)$. Let $P_1$ be the join of $4$-gon with tetrahedron and $P_2$ the simplicial polytope with $16$ facets. The polytope $P_1$ has $8$ facets and $26$ codimension $2$-faces while the polytope $P_1$ has $16$ facets and $48$ codimension $2$-faces. Performing the connected sum between these two polytopes we obtain a polytope with $22=16+8-2$ facets and $68=48+26-6$ codimension $2$-faces. Hence $(22, 68) \in \mathcal{E}^6$.
Perfoming the same computation with the polytope $P_1$ and the polytope with $11$ facets and $35$ codimensional $2$-faces we obtain $(25, 77) \in \mathcal{E}^6$. Finally the same operation between the polytopes $P_1$ and the simplicial polytope with $24$ facets gives $30=24+8-2$, $92=72+26-6$ which ends the proof that 
$(30, 92) \in \mathcal{E}^6$. This can also be proved by truncating a simple vertex from the $6$-polytope with $25$ vertices and $77$ edges. 
%
%
%
%
%
\qed

\section{Characterization of all $(f_0, f_1)$-pairs of $7$-polytopes with excess larger than $11$}\label{sect:four}
Starting from a $d$-dimension polytope for $d=3, 4, 5, 6$ we realize that the more difficult cases to check in order to characterize the $f$-vectors or at least the first two entries are related to the minimal values of the excess degree, 0, $d-2$, $d-1$ or $d$. For example the result becomes nicer if we decide to characterize the $(f_0, f_1)$-pairs for polytopes with large excess degree.

Consider the $4$-dimensional convex polytopes. It has been proved in \cite{BGrunbaum} that $$\mathcal{E}^4=\left\{(f_0, f_1): 2f_0\leq f_1\leq \binom{f_0}{2}\right\}\setminus\Bigg\{(6, 12), (7, 14), (8, 17), (10,20)\Bigg\}$$
None of the polytopes with
 $(f_0, f_1)$ pair $(6, 12), (7, 14), (8, 17), (10,20)$ has excess degree greater than $3\times 4-10=2$. 

Let  $\mathcal{E}_{>3d-10}^d$ be the set of $d$-polytopes whose excess degree is larger than $3d-10$. For $d=4$, the set  $\mathcal{E}_{>2}^4$ of $4$-polytopes whose excess degree is larger than $2$ is given by: $$\mathcal{E}_{>2}^4=\left\{(f_0, f_1): 1+2f_0< f_1\leq \binom{f_0}{2}\right\}.$$ In the same way for $d=5, 6$ we obtain:
$$\mathcal{E}_{>5}^5=\left\{(f_0, f_1): \frac{5}{2}+\frac{5}{2}f_0< f_1\leq \binom{f_0}{2}\right\},$$
and
$$\mathcal{E}_{>8}^6=\left\{(f_0, f_1): 4+3f_0< f_1\leq \binom{f_0}{2}\right\}.$$

Our goal in this section is to generalize the previous results for $d$-polytopes with excess degree great than $3d-10$. We can remark that for two integers $k$, $k'$ is $k<k'$ then $\mathcal{E}_{>k}^d\subset \mathcal{E}_{>k'}^d$ and then to find a complete characterization we just need to minimize $k$ in $\mathcal{E}_{>k}^d$

We now investigate the case of $7$-polytopes and conjecture an answer for a general $d$-polytope. We recall that if $P$ is a $d$-polytope having a simple vertex $v$ and $Q$ a $d$-polytope obtained from $P$  by truncating the vertex $v$ then $$f_0(Q)=f_0(P)+d-1\mbox{ and }  f_1(Q)=f_0(P)+\binom{d}{2}.$$ If for a $d$-polytope $P$ we have $f_1(P)\leq \frac{d+1}{2}f_0(P)$ then $P$ has at least one simple vertex.

Let $d\geq 4$ be an integer. In the remaining of this work we set $\epsilon_d(v)=2q-dp$ for any $v=(p, q)$ with $p\geq d+1$ and $\frac{d}{2}p\leq q\leq \binom{p}{2}$. This is nothing but the excess degree of a $d$-polytope having $v$ as $(f_0, f_1)$-pair if it exists.

\begin{theorem} \label{theo:dimsept} Let $\mathcal{E}^7$ be the set of $(f_0, f_1)$-pairs of $7$-polytopes. For $v=(p, q)$ such that $p\geq 8$ and $\frac{7}{2}p\leq q\leq \binom{p}{2}$, if $v\notin \mathcal{E}^7$ then $\epsilon_7(v)\leq 4\times7-10=11$. In other words the set of $(f_0, f_1)$-vector pairs for $7$-polytopes with 
excess strictly larger than $11$ is given by $$\mathcal{E}_{>11}^7=\Bigg\{(f_0, f_1): \frac{7}{2}f_0+\frac{11}{2}<f_1\leq  \binom{f_0}{2}\Bigg\}.$$
\end{theorem}

\proof
From the previous section we had
\begin{eqnarray}
\mathcal{E}^6=\left\{(f_0, f_1): 3f_0\leq f_1\leq \binom{f_0}{2}\right\}\setminus \Bigg(\left\{\Big(f_0, 3f_0+1\Big): f_0\geq 7\right\}\cup&&\cr \Big\{(8, 24); (9, 27); (9, 29); (10, 30); (10, 32); (10, 34); (11, 33);  (12, 38); (12, 39); (13, 39); (14, 42);&&\cr  (14, 44); (15, 47); (18, 54); (19, 57);(17, 53); (20, 62)\Big\}\Bigg)\cr
\end{eqnarray}

A pyramid over the $6$-polytopes gives:
\begin{eqnarray}
\left\{(f_0, f_1): 4f_0-4\leq f_1\leq \binom{f_0}{2}\right\}\setminus \Bigg(\left\{\Big(f_0+1, 4f_0+1\Big): f_0\geq 7\right\}\cup&&\cr \Big\{(9, 32); (10, 36); (10, 38); (11, 40); (11, 42); (11, 44); (12, 44);  (13, 50); (13, 51); (14, 52); (15, 56);&&\cr  (15, 58); (16, 62); (18, 70); (19, 72); (20, 76); (21, 82)\Big\}\Bigg)\subset \mathcal{E}^7.
\end{eqnarray}
A direct computation shows that $\epsilon((f_0+1, 4f_0+1))>11$ if and only if $f_0>17$. Assume that $f_0>17$ and let us prove that $(f_0-6, 3f_0-14)\in \mathcal{E}^6$. We have $\epsilon_6((f_0-6, 4f_0-20))=8$ and if $(f_0-6, 3f_0-14)\notin \mathcal{E}^6$
then $(f_0-6, 3f_0-14)=(10, 34)$, as $(10, 34)$ is the only vector not in $\mathcal{E}^6$ with excess equal to $8$. Therefore we get $f_0=16$ which is a contradiction. In conclusion for $f_0>17$ there is a $6$-polytope $P$ with $(f_0, f_1)$-pair $(f_0-6, 3f_0-14)$; and a pyramide over $P$ give a $7$-polytope $Q$ having $(f_0, f_1)$-vector which is equal to $(f_0-5, 4f_0-20)$. As $4(f_0-5)<( 4f_0-20)+1$ the polytope $Q$ has a simple vertex whose truncation gives a $7$-polytope having $(f_0, f_1)$-pair equals $(f_0+1, 4f_0+1)$. We can conclude that all the $7$-polytopes with excess greater than $11$ and with $(f_0, f_1)$-pairs in $\Big\{\big(f_0+1, 4f_0+1\big): f_0\geq 7\Big\}$ exist.

Let us focus on the set 
\begin{eqnarray}
L=\Big\{(9, 32); (10, 36); (10, 38); (11, 40); (11, 42); (11, 44); (12, 44);  (13, 50); (13, 51); (14, 52);\cr  (15, 56); (15, 58); (16, 62); (18, 70); (19, 72); (20, 76); (21, 82)\Big\}.
\end{eqnarray}
The only vectors $v=(p, q)\in L$ with $\epsilon_7(v)>11$ are $$v=(p, q)=(16, 62); (18, 70); (20, 76); (21, 82); (23, 90); (26, 102); (31, 123).$$
For $v=(p, q)=(16, 62); (20, 76); (21, 82); (26, 102); (31, 123)$, we compute $v'=(p-8, q-p-20)=(8, 26); (12, 36); (13, 41); (18, 56); (23, 72) \in \mathcal{E}^6$. Then their exist $6$-polytopes $P_{v'}$ whose $(f_0, f_1)$-pairs are equal to $v'$. A pyramid over them give $7$-polytopes having $(f_0, f_1)$-pairs equal to $(p-7, q-28)=(9, 34); (13, 48); (14, 54); (19, 74); (24, 95)$. In each case we observe that $q-28<4(p-7)$ which means that each of them has a simple vertex whose truncation give $7$-polytopes with $(f_0, f_1)$-pairs equal to $(p-1, q-7)$. As truncations of simple vertices generate simplex facets then pyramids on these give the result.

Let us now concentrate on the cases $v=(18, 70)$. We first consider $v=(18, 70)$. There is a $6$-polytope $R$ with  $(f_0, f_1)=(10, 35)$. A pyramid over $R$ gives a $7$-polytope $R'$ having  $(f_0, f_1)$-vector equal to $(11, 42)$. As $42<4\times 11$ then $R'$ has a simple vertex whose truncation gives a $7$-polytopes $R''$ with  $(f_0(R''), f_1(R''))=(17, 63)$. The truncation of a simple vertex in $R''$ with generate a simplex facet $F$ and a pyramid other $F$ gives a $7$-polytope with $(f_0, f_1)$-vector equal to $(18, 70)$. 
%

We now turn to the pair $v=(f_0, f_1)$ with $f_0\geq 8$ and $f_1\in ]\frac{7}{2}f_0, 4f_0 + 1[$. The condition $\epsilon_7(v)>11$ implies that $f_1\geq \frac{11}{2}+\frac{7}{2}f_0$ and then we need to discuss two cases: $\frac{11}{2}+\frac{7}{2}f_0>4f_0-4$ and $\frac{11}{2}+\frac{7}{2}f_0<4f_0-4$.

If $\frac{11}{2}+\frac{7}{2}f_0>4f_0-4$ then there is nothing else to prove as we end up in the pyramid case. Suppose that $\frac{11}{2}+\frac{7}{2}f_0<4f_0-4$ i.e. $f_0>19$ and set for $k\geq 8$, $X^7_k=\{(k, f_1); \frac{11}{2}+\frac{7}{2}f_0< f_1 < 4f_0-4\}$. We can prove by truncation that if $X^7_k\subset \mathcal{E}_{>11}^7$, then $X^7_{k+6}\subset \mathcal{E}_{>11}^7$. To prove that each vector $(f_0, f_1)$ satisfying this condition defines a $7$-polytope it is sufficient to show that $X^7_k\subset  \mathcal{E}_{>11}^7$ for $k=8,\cdots, 13$. But $k\leq 19$ implies that
$\frac{11}{2}+\frac{7}{2}k>4k-4$ which have already been solved.

Finally we conclude that all the pairs $(p, q)$ with  $p\geq 8$, $\epsilon_7(v)>11$ and $\frac{7}{2}p\leq q\leq \binom{p}{2}$, characterize $7$-polytopes. In other words
the set of $(f_0, f_1)$-vectors pair for $7$-polytopes with 
excess strictly larger than $11$ is given by $$\mathcal{E}_{>11}^7=\Bigg\{(f_0, f_1): \frac{7}{2}f_0+\frac{11}{2}<f_1\leq  \binom{f_0}{2}\Bigg\}.$$

\qed
%
%

From all the previous results we can make the following conjecture.

%
%

\begin{conj} \label{conj:dimarb} Let $d\geq 4$ be an integer and $\mathcal{E}^d$ be the set of $(f_0, f_1)$-pairs of $d$-polytopes. For $v=(p, q)$ such that $p\geq d+1$ and $\frac{d}{2}p\leq q\leq \binom{p}{2}$, if $v\notin \mathcal{E}^d$ then $2q-dp\leq 4d-10$. In other words the set of $(f_0, f_1)$-pairs for $d$-polytopes; $d\geq 4$ with 
excess strictly larger than $3d-10$ is given by $$\mathcal{E}_{>3d-10}^d=\Bigg\{(f_0, f_1): \frac{d}{2}f_0+\frac{3d-10}{2}<f_1\leq  \binom{f_0}{2}\Bigg\}.$$
\end{conj}
%
%
%
%

\section*{Acknowledgments}
The research leading to these results has received funding from the European Research Council under the European Union's Seventh Framework Programme ERC Grant agreement ERC StG 716424 - CASe and the Israel Science Foundation under ISF Grant 1050/16. RCA research at Max-Planck Institute is supported by the Alexander von Humboldt foundation.

\vspace{0.5cm}


\end{document}